\documentclass[a4paper,12pt,leqno]{amsart}

\usepackage{amssymb,hyperref}
\usepackage[latin1]{inputenc}
\swapnumbers
\theoremstyle{plain}

\newtheorem{theo}{Theorem}[section]
\newtheorem{prop}[theo]{Proposition}
\newtheorem{cor}[theo]{Corollary}

\theoremstyle{definition}

\theoremstyle{remark}
\newtheorem{rem}[theo]{Remark}

\numberwithin{equation}{section}

\newcommand{\thismonth}{\ifcase\month\or
  January\or February\or March\or April\or May\or June\or July\or
  August\or September\or October\or November\or December\fi
  \space\number\year}
\DeclareMathAlphabet{\mathrmsl}{OT1}{cmr}{m}{sl}

\newcommand{\oper}[3][n]{\newcommand{#2}{\mathop{\mathrm{#3}}%
\ifx n#1\nolimits\else\limits\fi} }
\newcommand{\rsoper}[3][n]{\newcommand{#2}{\mathop{\mathrmsl{#3}}%
\ifx n#1\nolimits\else\limits\fi} }
\newcommand{\vol}{\operatorname{vol}}
\newcommand{\Scal}{\operatorname{Scal}}

\newcommand{\csum}{\mathfrak{S}}
\newcommand{\asum}{\mathfrak{A}}
\newcommand{\tv}{\mathcal{T}}

\newcommand{\II}{\mathbb{I}}

\renewcommand{\exp}{\operatorname{e}}
\newcommand{\gb}{\overline{g}}
\newcommand{\Id}{\operatorname{Id}}
\newcommand{\ti}{\theta^1}
\newcommand{\vti}{\vartheta^1}
\newcommand{\vt}{\vartheta^0}
\newcommand{\vtb}{\vartheta^{\bar{0}}}
\newcommand{\tib}{\theta^{\bar{1}}}

\newcommand{\fti}{\hat{\theta}^1}
\newcommand{\ftib}{\hat{\theta}^{\bar{1}}}
\newcommand{\taui}{\tau^1}
\newcommand{\tauib}{\tau^{\bar{1}}}

\newcommand{\proofof}[1]{\end{#1}\begin{proof}}

\newcounter{mnotecount}[section]
\renewcommand{\themnotecount}{\thesection.\arabic{mnotecount}}
\newcommand{\mnote}[1]
{\protect{\stepcounter{mnotecount}}$^{\mbox{\footnotesize  $
      \bullet$\themnotecount}}$ \marginpar{\raggedright\tiny\em
    $\!\!\!\!\!\!\,\bullet$\themnotecount: #1} }

\begin{document}

\title[Renormalized volume for ACHE manifolds]{A remark on renormalized volume\\ 
and Euler characteristic for ACHE $4$-manifolds}
\author{Marc Herzlich}
\address{Institut de Math\'ematiques et Mod\'elisation de Montpellier\\ 
UMR 5149 CNRS -- Université Mont­pellier~II\\ France}
\email{herzlich@math.univ-montp2.fr}

\begin{abstract}
This note computes the ``renormalized volume'' and a 
renormalized Gauss-Bonnet-Chern formula for asymptotically complex hyperbolic 
Einstein (the so-called {\sc ache}) $4$-manifolds.
\end{abstract}
\thanks{The author is a member of the {\sc edge} Research Training 
Network {\sc hprn-ct-2000-00101} of the European Union  and is supported 
in part by an {\sc aci} program of the French Ministry of Research.}
\maketitle
\section{Introduction.}

Asymptotically symmetric Einstein metrics exhibit many interesting 
phenomena~\cite{biquard-symmetric,graham-lee}. They were especially studied in the
asymptotically {\sl real hyperbolic} (or {\sc ahe}) case, which enjoys fruitful 
relationships with 
physics through the {\sc ads-cft} correspondence. They are also a useful tool
in the study of conformal geometry in establishing links between the
conformal geometry of a compact $(n-1)$-dimensional manifold (usually called 
the {\sl boundary at infinity}) and a complete Einstein $n$-dimensional manifold
(the {\sc ahe} manifold). In this setting, an intriguing 
invariant, called {\sl renormalized volume}, has been defined by C.~R.~Graham 
\cite{graham-volume}, after works by physicists such as Henningson and 
Skenderis~\cite{henningson-skenderis}. 

In even dimensions $n$, the renormalized 
volume is an invariant of the Einstein metric only. If $n=4$, its role in the 
formula for the Euler characteristic of the Einstein manifold has been moreover
pointed out by M.~T.~Anderson \cite{anderson-volume}, with applications in the 
study of the moduli space of Einstein asymptotically real hyperbolic metrics 
\cite{anderson-existence}. This formula is called ``the renormalized 
Gauss-Bonnet-Chern formula'': although the Einstein manifold is non-compact,
all divergent terms in the integrals of the formula are shown to cancel, whereas
renormalized volume appears as a finite limit contribution. 

In odd dimensions
$n$, the renormalized volume is not an invariant of the Einstein metric
only but rather depends on a choice of a representative metric on the boundary 
at infinity in the conformal class. This makes it no less interesting, as it gives
rise to the so-called {\sl conformal anomaly} phenomenon: the difference between
the renormalized volumes of two different choices of metric singles out
a local differential operator on the boundary with nice properties 
\cite{graham-volume}.

\medskip

The goal of this short note is to point out analogous results in the case of
Einstein asymptotically {\sl complex} hyperbolic (or {\sc ache}) manifolds of
dimension $4$,
where the boundary at infinity is now a strictly pseudoconvex $3$-manifold,
with the hope that such an object would be interesting for the study of
$3$-dimensional CR geometry.

Unfortunately, the situation is less pleasant than in the real case, as the 
renormalized volume is never an invariant of the the complete Einstein metric and 
always depends on the choice of some contact form (or in the 
usual language of CR geometry: a \emph{pseudo-hermitian structure}) compatible
with the CR-structure at infinity. This situation is reminiscent from that of
Einstein asymptotically {\sl real} hyperbolic manifolds of odd-dimensions 
({\sl i.e.} boundary at infinity of {\sl even} dimension), this
should come as no surprise as it is well-known that CR geometry enjoys lots
of analogies with even-dimensional conformal geometry.

However, the fact that in the {\sc ache} context the bulk Einstein manifold is 
even-dimensional brings some nice features. It turns out that adding some 
well-chosen local quantity at infinity can yield an invariant of the Einstein
metric only. As expected, a renormalized Gauss-Bonnet characteristic formula 
can be obtained as well.

Our main results then read as follows (further notations and definitions are given
in the next section):

\medskip

\begin{theo}\label{theo1}
Let $(M,g)$ be a $4$-dimensional Einstein asymptotically complex hyperbolic 
\upn{(}{\sc ache}\upn{)} manifold, with boundary at infinity a compact strictly 
pseudoconvex CR 
$3$-dimensional manifold $X$ with contact distribution $H$ and almost-complex
structure $J$. Then there exists for any choice of compatible contact form 
$\eta$ on $X$ an invariant $V$ of the pair $(g,\eta)$ called \emph{renormalized
volume} of $g$ relative to $\eta$.
\end{theo}

\begin{theo}\label{theo2}
Under the same assumptions, if moreover $R$ and $\tau$ are the curvature and 
torsion of the Webster-Tanaka connection on $(X,H,J,\eta)$, then
\[ \mathcal{V} = \frac{3}{2}\, V - \int_X \left(\frac{R^2}{16} - 
\frac{5}{2} |\tau|^2 \right) \]
is an invariant of the metric $g$ only, and
\[ \chi(M) = \frac{1}{8\pi^2}\int_M \left( |W^g|^2 - \frac{1}{24}(\Scal^g)^2 \right)
+ \frac{1}{4\pi^2}\, \mathcal{V} .\]
\end{theo}

\medskip

As the model case of the complex hyperbolic plane shows, the appearance
of the integral factor on the boundary seems unavoidable; see section~2 for further 
details. This shows than, rather than giving rise to a global invariant, the 
renormalized volume gives birth to a {\sl conformal anomaly}, {\sl i.e.} a formula 
relating the renormalized volume for some choice of pseudo-hermitian structure at 
infinity to its expression for some other choice at infinity, through a {\sl local}
differential expression. Namely, if we let $V(\eta)$ be the renormalized volume for 
a choice of contact form at infinity,

\begin{cor} For each contact form $\eta$, there is a differential operator $P_{\eta}$
on $X$ such that, for any function $f$ on $X$ which never vanishes,
\[ V(f\eta) - V(\eta) \ = \ \int_X  P_{\eta}(f) \,\eta\wedge d\eta ,\]
\end{cor}

In the {\sl real hyperbolic} case \cite{graham-volume}, the conformal anomaly is given 
by differential operators with nice invariance properties. Our result in the {\sl complex 
hyperbolic case} strongly suggests that it should be interesting to study the operator 
arising from the variation of
\[ T(\eta) \ = \ \frac{R^2}{16} -  \frac{5}{2} |\tau|^2 \] 
under deformations of the contact form in the same contact structure. Explicit derivations
of the variations of Tanaka-Webster curvature $R$ and torsion $\tau$ are given as an 
Appendix to this note; further study will be deferred to a future work.

Note moreover that N. Seshadri has given in \cite{seshadri-vol} another version
of the renormalized volume that covers all dimensions but for K\"ahler-Einstein 
metrics only rather than ACHE.

\bigskip

\section{Definitions and notations.}

\medskip

Let $(X^3,H,J_0)$ be a strictly pseudo-convex 
$3$-dimensional CR manifold, {\sl i.e.} a contact 
manifold with contact distribution $H$ and almost complex structure $J_0$ on $H$. If 
$\eta$ in any choice of compatible contact form, an associated metric $\gamma$ may be 
defined on $H$ by $\gamma = d\eta(\cdot, J_0\cdot)$. One gets from it a Reeb field 
$\xi$ and a (Tanaka-Webster) connection $\nabla$ whose torsion in the direction of $\xi$ 
is $\tau = \nabla^{TW}_{\xi}\cdot - [\xi, \cdot]$.

\medskip

Let $M$ be a $4$-manifold such that 
the complement of some compact set is diffeomorphic to $[r_0,+\infty[\times X$.
We consider first the metric $g_0 = dr^2 + \exp^{2r}\eta^2 + \exp^r\gamma$ on 
$]r_0,+\infty[\times X$ and let $C^{\infty}_{\delta}$ be the space of smooth functions 
on $M$ such that $\exp^{\delta r}\nabla^k f$ is bounded for any $k$. 
Any metric $g$ on $M$ such that $g - (dr^2 + \exp^{2r}\eta^2 + \exp^r\gamma)$ belongs
to $C^{\infty}_{\delta}$ for some $\delta >0$ will be called an {\sl asymptotically 
complex hyperbolic} metric. Moreover, $(M,g)$ is said to be {\sc ache} if $g$ is an 
Einstein metric.

A lot of such metrics arise on pseudoconvex domains in $\mathbb{C}^2$ (and are 
K\"ahler-Einstein in this case \cite{cheng-yau}) whereas another important family 
was constructed by O.~Biquard in \cite{biquard-symmetric}. The Biquard metrics
are especially interesting in the case the boundary at infinity $X$ is endowed with a 
non-embeddable CR structure, as they provide a substitute for the non-existing
K\"ahler-Einstein metric.
 
\medskip

In \cite{ob-mh}, the author and O.~Biquard carefully studied the asymptotic behaviour of 
{\sc ache} metrics, and precise asymptotic expansions were obtained. In all that
follows, we consider an {\sc ache} metric $g$ on a neighbourhood of infinity 
$]r_0,+\infty[\times X$ in $M$. If a contact form $\eta$ is given, there exists a
canonical Tanaka-Webster connection $\nabla$ on $X$. For any tensor field $D$ on $X$,
$D_{a,bc...}^{d...}$ ($a,b,c,d,... \in \{1,\bar{1}\})$ will denote the components of 
$D$ (and subsequent Tanaka-Webster derivatives, separated by a comma from the original
components) in a local orthonormal coframe $(\ti,\tib)$, {\it i.e.} such that
$d\eta = i\ti\wedge\tib$. For instance we shall use expressions such 
as $\tau^{a}_{b,c...}$ for the (derivatives of the) torsion $\tau$ of $\nabla$ and 
$R_{,ab...}$ for its curvature. We also denote $\vt:=\exp^{-r}dr+i\eta$, and 
$\vti:=\ti-\phi\lrcorner\ti$. Last, in any power series expansion $\sum \phi_k(x)
\exp{kr}$, the $k$-the term $\phi_k$ (seen as a function on $X$) will be called 
{\sl formally determined} if it can be computed with the knowledge of a finite jet
of the CR structure at $x\in X$ only. The most interesting feature of {\sc ache} metrics
(and K\"ahler-Einstein metrics as well) is that they are not entirely formally
determined. The results in \cite{ob-mh} are summarized in the three  following 
statements:

\begin{theo}[\cite{ob-mh}]\label{theo-expansion} 
There exists on $]r_0,+\infty[\times X$ an integrable complex structure $J$ given by a 
\upn{(}not necessarily convergent\upn{)} formal series, entirely determined formally from 
data at infinity. The first terms in its expansion is given by
\[ J = J_0 -2 \exp^{-r}\tau + \exp^{-2r}(2|\tau|^2 - J_0\nabla_\xi\tau) + o(\exp^{-2r}), \]
or equivalently by a map $\phi = -i\exp^{-r}\tau + 
\frac{1}{2}\exp^{-2r}\nabla_{\xi}\tau + o(\exp^{-2r})$ from $T^{0,1}_{J_0}$ to 
$T^{1,0}_{J_0}$.
\end{theo}

\begin{theo}[\cite{ob-mh}]\label{theo-expan2}  
There is  on $]r_0,+\infty[\times X$ a \upn{(}formal series\upn{)} K\"ahler-Einstein 
metric $\gb$. The K\"ahler form $\omega$ of $\gb$ is formally determined up to order 
$2$ as follows
\begin{equation*}\begin{split}
\omega \ = \ & \exp^r \left( dr\land\eta + d\eta\right)
- \frac{R}{2}\, d\eta \\ 
& + \frac{4}{3}\left(\frac{i}{8}R_{,\bar 1}\vt\land\tib 
- \frac{i}{8}R_{,1}\vtb\land\ti -\frac{1}{2}\taui_{\bar 1,1}\vt\land\tib
-\frac{1}{2}\tauib_{1,\bar 1}\vtb\land\ti\right) - \frac{\Delta R}{2}\exp^{-r}d\eta\\
& - \frac{2}{3}\left(\frac{R^2}{8}-|\tau|^2-\frac{\Delta R}{6}
 +\frac{2i}{3}(\taui_{\bar{1},11}-\tauib_{1,\bar{1}\bar{1}})\right)\exp^{-r}dr\land\eta\\
& + \frac{2}{3}\left(\frac{R^2}{8}-|\tau|^2 +\frac{\Delta R}{12}
 -\frac{i}{3}(\taui_{\bar{1},11}-\tauib_{1,\bar{1}\bar{1}})\right)\exp^{-r}d\eta + 
o(\exp^{-2r}).
\end{split}\end{equation*}
Moreover, if $g$ is an {\sc ache} metric with the same boundary at infinity, then there 
exists an anti-$J_0$-invariant symmetric bilinear form $k$ on 
$H$ and a unique diffeomorphism $\psi$ asymptotic to identity at infinity such that 
$\psi^*g = \gb + k \exp^{-r} + o(\exp^{-2r})$.
\end{theo}

\begin{cor}[\cite{ob-mh}]\label{cor-expan} 
The K\"ahler metric $\gb$ is explicitly given by  
\begin{equation*}\begin{split} 
\gb \ =\ & (dr^2 + \exp^{2r}\eta^2 + \exp^r \gamma)  - \frac{R}{2}\gamma + 
2\gamma(J_0\tau\cdot ,\cdot ) + \frac{1}{6}(R_{,1}\ti\circ\vtb + 
R_{,\bar 1}\tib\circ\vt)\\
& +\frac{2i}{3}(\taui_{\bar 1,1}\vt\circ\tib - \tauib_{1,\bar 1}\vtb\circ\ti) 
- \exp^{-r}R\gamma(J_0\tau\cdot ,\cdot)
- \exp^{-r}\gamma (\nabla_{\xi}\tau(\cdot),\cdot )\\ 
& - \frac{2}{3}\left(\frac{R^2}{8}-|\tau|^2-\frac{\Delta R}{6}
 +\frac{2i}{3}(\taui_{\bar{1},11}-\tauib_{1,\bar{1}\bar{1}})\right)
\exp^{-2r}(dr^2+\exp^{2r}\eta^2)\\ 
& + \frac{2}{3}\left(\frac{R^2}{8}-|\tau|^2 +\frac{\Delta R}{12}
-\frac{i}{3}(\taui_{\bar{1},11}-\tauib_{1,\bar{1}\bar{1}})\right)\exp^{-r}\gamma + 
o(\exp^{-2r}),
\end{split}\end{equation*}
where $\alpha\circ\beta = \alpha\otimes\beta + \beta\otimes\alpha$ is the 
symmetrized product of forms. 
\end{cor}

The main conclusion of these facts is the following : given any
ACHE metric $g$ and any choice of pseudo-hermitian structure at infinity realizing
the CR structure, there is a unique diffeomorphism $\psi$ asymptotic
to the identity on $X$ such that, up to strictly lower order terms, $\psi^*g$ can be 
written as the sum of a formally determined K\"ahler-Einstein metric and a 
{\sl formally undetermined} term of order $2$ (notice that $k\exp^{-r}$ decays
like $\exp^{-2r}$). From now on, we will forget the 
diffeomorphism $\psi$ and, if the metric $g$ is written this way
in such coordinates, we will say that it is ``in the K\"ahler gauge associated 
to the choice of pseudo-hermitian structure at infinity''. 

Using the results of \cite{ob-mh} that we have just recalled, we can now define the
renormalized volume:

\begin{prop} \label{prop-vol}
Let $g$ be an ACHE metric on $M$, written in a K\"ahler gauge associated to some choice
of pseudo-hermitian structure at infinity. Then the volume of large coordinate balls 
$B(r)$ of radius $r$ \upn{(}complement of $]r,+\infty[\times X$ in $M$\upn{)} 
has an asymptotic expansion: $\vol_g B(r) = \pi^2 \exp^{2r} +\, v_1 \exp^r 
+\, V +\, o(1)$. the number $V$ is the \emph{renormalized volume} of the metric $g$ 
associated to the choice of pseudo-hermitian structure at infinity
\end{prop}

\begin{proof}  
To check the proposition, just notice that the volume form of $g$ only differs from
that of $\gb$ at order $\frac{5}{2}$ since $k$ is trace-free, and, in the volume form 
of the K\"ahler form $\omega$, order $\frac{3}{2}$ terms do not exist whereas order 
$2$ terms are of zero integral from the CR Stokes' formula \cite{cheng-lee}: whenever
$\alpha=\alpha_1\ti$ is a $(1,0)$-form on $X$ (given in a local orthonormal coframe), 
one has
\[ d\alpha = \alpha_{1,\bar 1}\ti\wedge\tib + \alpha_{1,0}\eta\wedge\ti + 
\alpha_1\tau_{\bar 1}^1\eta\wedge\tib \]
(recall $\alpha_{\cdot,\cdot}$ denotes the components of $\nabla\alpha$ in the
local coframe), and 
\[ \int_X \alpha_{1,\bar 1}\eta\wedge\ti\wedge\tib = \int_X (d\alpha)\wedge \eta =
 - \int_X \alpha\wedge d\eta = 0 ,\]
since $X$ is closed and $d\eta=i\ti\wedge\tib$. This achieves the proof of Theorem
\ref{theo1}.\end{proof}
 
\medskip

\begin{rem}
In the asymptotically {\sl real} hyperbolic Einstein case ({\sc ahe}), the renormalized 
volume
is similarly defined \cite{graham-volume}, but with the help of a different gauge. 
It is proved in \cite{graham-volume} that it is always possible, given a metric 
$h^{\infty}$ in the conformal class of the boundary at infinity,
to find coordinates such that $g=dr^2 + h(r)$ on $]r_0,+\infty[\times X$ (this is
the {\sl geodesic gauge}, the function $r$ being the {\sl geodesic defining function} 
associated to a choice of metric $h^{\infty}$ in the conformal class at infinity).
The metric $h(r)$ has an expansion in powers of $\exp^r$ and $V$ is defined as above as 
the constant coefficient in the expansion of $\vol(B(r))$. It is is in itself an invariant 
of $g$. The reader might hence think that the ``misbehaviour'' of the renormalized volume 
in the {\sc ache} case comes from a bad choice of gauge. However, the standard metric of
$\mathbb{C}\mathbf{H}^2$ is both in the K\"ahler and geodesic gauges, and Theorem
\ref{theo1} yields that the boundary term in the renormalized Gauss-Bonnet formula
(which necessarily is an invariant of $g$) {\sl is different} than the renormalized
volume.
\end{rem}

\begin{rem} The most important fact to be noted in the previous Proposition is that
{\sl there is no term in} $\vol_g B(r)$ {\sl that is linear in} $r$. In the {\sc ahe} 
case \cite{graham-volume}, when the boundary at infinity is {\sl odd-dimensional},
an analogous phenomenon occurs: linear terms cancel in the asymptotic expansion of large 
balls. In the {\sc ahe} case again but when the boundary at infinity is 
{\sl even-dimensional} (the case that is considered to be the closest to the {\sc ache} case,
although dimensions of the boundaries at infinity differ), the situation is different: some
non-zero linear term appears in the expansion of the volume, with a coefficient related to 
the integral term in the Gauss-Bonnet-Chern formula for the boundary at infinity $X$
itself \cite{graham-volume}. 

If one believes in this analogy (between {\sl even}-dimensional conformal geometry and CR
geometry), one might then wonder why there is indeed no linear term in Proposition 
\ref{prop-vol}. However, reasoning by analogy again, one would assert from \cite{graham-volume}
that the coefficient of the linear term should be a multiple of the {\sl integral} 
$Q$-{\sl curvature} of $X$ \cite{fefferman-hirachi}, but it has been proved in 
\cite{fefferman-hirachi} that this integral always vanishes in $3$-dimensional CR 
geometry, thus the absence of any linear term, a phenomenon that might be purely
$3$-dimensional.\end{rem}

\bigskip

\section{Proof of Theorem \ref{theo2}.}

\medskip

We first choose a contact form (or pseudo-hermitian structure) at infinity 
({\sl i.e.} on $X$) realizing 
the CR structure and we put the ACHE metric in the associated K\"ahler gauge around
infinity. The basic element of the proof then is the Gauss-Bonnet-Chern formula for 
the Euler characteristic of the compact domain with boundary $B(r)$ 
delimited by what we shall call the 
{\sl coordinate sphere} $S(r)= \{r\}\times X$:
\begin{equation}\label{chi}\begin{split}
\chi(B(r)) \ = \ & \frac{1}{8\pi^2} \int_{B(r)} \left( |W|^2 - \frac{1}{24}\Scal^2\right)
\ + \ \frac{1}{96\pi^2}\int_{B(r)} \Scal^2 \\
& + \ \frac{1}{12\pi^2}
\int_{S(r)} \tv(\II\wedge\II\wedge\II) + 3\, \tv(\II\wedge \mathcal{R}) ,
\end{split}\end{equation}
where the $\wedge$ operation provides a $p+q$-form  with
values in $\otimes^{r+s}TM$ from a $p$-form with values in $\otimes^rTM$ and a $q$-form
with values in $\otimes^sTM$, and we have denoted by $\tv$ the contraction between
the volume form of $S(r)$ and elements of $\otimes^3TM$. Moreover, $\II$ is the shape 
operator of $S(r)$ in $(M,g)$ and $\mathcal{R}$ is the curvature ($2$-form with values 
in $2$-forms) of $(M,g)$, with $W$ and $\Scal$ denoting its Weyl and scalar curvature
(trace-free Ricci curvature is zero as $g$ is Einstein).
Notice also the difference in notation between the (scalar) curvature $R$ of the
$3$-dimensional CR manifold $X$ and the curvature tensor $\mathcal{R}$ of the 
$4$-dimensional Einstein manifold $M$.

\medskip

It is proven in \cite{ob-mh} that the integral involving $|W|^2 - \frac{1}{24}\Scal^2$
on $B(r)$ converges for an {\sc ache} metric when $r$ goes to infinity. Moreover, it 
is clear that both the scalar curvature integral (which is, up to 
a constant, $\vol B(r)$)  and the boundary integrals have an asymptotic expansion
in powers of $\exp^{-\frac{r}{2}}$ (there are no polynomial terms as they cancel
in the volume expansion, as noted above). Convergence of all the other terms implies
that divergent terms cancel pairwise, whereas the limit as $r$ goes to infinity
of 
\[ \chi(M) - \int_{B(r)} \left( |W|^2 - \frac{1}{24}\Scal^2\right) - 
\frac{3}{8\pi^2}V \]
is given by the constant terms in the asymptotic expansion of the boundary integrals.
Our task then reduces to a careful computation of these terms. For this, the following 
facts will be useful:

\smallskip

{\flushleft {\bf Fact 1}}. It is proven in \cite{ob-mh} that replacing $g$ by $\gb$ in 
the boundary integrals introduces terms that are $o(\exp^{-2r})$
only, hence do not contribute in the limit as the volume form of each sphere
is $O(\exp^{2r})$ at most. Hence all computations can be done using the formal
K\"ahler-Einstein metric $\gb$ rather than the ACHE metric $g$.

\smallskip

{\flushleft {\bf Fact 2}}. More precisely, the highest-order term in the
expansion of the volume form of $S(r)$ is $\exp^{2r}\eta\wedge d\eta$
(where $\eta$ is the contact form underlying the chosen pseudo-hermitian structure
chosen at infinity). 
Hence we will only need to track the order $\exp^{-2r}$ terms in the proof
below. Every asymptotic expansion we will meet in the course of the computations 
is of the following type:
\[ A_0 \exp^{2r} \ + \ A_1\exp^{-r} \ + \ A_{\frac{3}{2}} \exp^{-\frac{3}{2}r}
\ + \ A_2 \exp^{-2r} \ + \ o(\exp^{-2r}). \]
As a result, order $2$ terms may arise during the computation
only when putting together an order $2$ term with
order $0$ terms or two order $1$ terms with order $0$ terms. Order $\frac{3}{2}$ terms 
can hence be forgotten during the whole computation, {\sl unless} when some 
differentiation is involved, as doing so along directions in $X$ raises the order 
possibly by a factor $\frac{1}{2}$.

\smallskip

{\flushleft {\bf Fact 3}}. Our final computation involves integration along $X$, 
hence each exact term can be forgotten. Using the CR Stokes' formula already described in
the proof of Proposition \ref{prop-vol}, see also \cite{cheng-lee}, this means that 
every term in $R_{,1\bar 1}$, $R_{,\bar 1 1}$, $\Delta R=R_{,1\bar 1} + R_{,\bar 1 1}$, 
$\tau^1_{\bar 1,11}$ or $\tau^{\bar 1}_{1,\bar 1 \bar 1}$ drops out. In what follows, 
occurrence of such a term will be denoted by $\mathcal{O}$.

\smallskip

{\flushleft {\bf Fact 4}}. From \cite{ob-mh} again, the curvature tensor $\mathcal{R}$ of
the formal K\"ahler-Einstein metric (seen as a $2$-form with values in $2$-forms) is, up 
to order~$2$, given by the sum of the model curvature tensor ({\sl i.e.} that has the 
same expression w.r.t. $\gb,J$ as the constant holomorphic sectional curvature has
w.r.t. $g^{\mathbb{C}\mathbf{H}^2},J_0$) and of an order $2$ term, called $W_2^-$ and
controlled by the Cartan tensor of the CR-structure at infinity. Said shortly, one
writes: $\mathcal{R} = \mathcal{R}_0 + W_2^- \exp^{-2r}$.

\medskip

From now on, the task can be divided into three steps: computation of the outer unit
normal and intrinsic volume form of $S(r)$, computation of the shape operator
(the only step that involves differentiation) and estimation of the order $2$ terms 
in $\tv(\II\wedge\II\wedge\II)$ and $\tv(\II\wedge \mathcal{R})$. As the computations
involved are rather long, we shall give here the main intermediate results only, 
indicating at each stage which are the key steps and facts that lead to them.

\medskip

From the explicit expansion of $\gb$, we can get immediately the {\sl outer unit normal} 
of $S(r)$:
\begin{equation}\label{normale} 
\nu(r) \ = \ \left( 1 + \frac{1}{3}\exp^{-2r}\left( \frac{R^2}{8} - |\tau|^2 + 
\mathcal{O}\right) \right) \partial_r \ + \ \nu^{T} \ + \ o(\exp^{-2r}) 
\end{equation}
where $\nu^T$ is an order $\frac{3}{2}$ term tangent to $X$, involving linearly 
$R_{,1}$, $R_{,\bar 1}$, $\tau^1_{\bar 1,1}$ and $\tau^{\bar 1}_{1,\bar 1}$. It
will be proved below that it is not necessary to detail further the expression of this
term.

The {\sl volume form} $\varpi$ of $S(r)$ is then (up to forgotten order $\frac{3}{2}$ 
terms, see Fact 2 above):
\begin{equation}\label{volform}\begin{split} 
\varpi  & \ = \ \frac{1}{2} \, \omega^2( \nu(r), ., ., .) \\
       & = \ \exp^{2r}(1 + \exp^{-r}\varpi_1 + \exp^{-2r}\varpi_2)\, \eta\wedge d\eta
 \ + \ o(\exp^{-2r}) \\
        & \ = \ \exp^{2r}\left( 1 - \frac{R}{2}\exp^{-r} + \frac{1}{3}\exp^{-2r}
\left(\frac{R^2}{8} - |\tau|^2 + \mathcal{O} \right) \right)\, \eta\wedge d\eta \ + 
\ o(\exp^{-2r}).\end{split}\end{equation}

The {\sl shape operator} $\II$ is obtained by taking the 
extrinsic covariant derivative of the unit outer normal $\bar{\nabla}\nu(r)$
(where $\bar{\nabla}$ here denotes the Levi-Civita connection of $\gb$. As $\nu^T$ is 
an order $\frac{3}{2}$ term, only its derivatives in the direction of $H$ might 
contribute to order $2$ terms in $\II$, but it is an easy task to convince oneself 
that these would add only terms linear in $R_{,1\bar 1}$, $R_{,\bar{1}1}$, 
$\tau^1_{\bar 1,11}$ and $\tau^{\bar 1}_{1,\bar 1\bar 1}$, hence of vanishing integral
from Fact 3. As a result, they can be forgotten.

It remains to compute the derivative of the radial term in $\nu(r)$, seen first as a
bilinear symmetric form. Keeping only symmetric terms in the usual $6$-term formula 
for the covariant derivative, the only significant term is
\[ \frac{1}{2}\left( 1 + \frac{1}{3}\exp^{-2r}\left( \frac{R^2}{8} - |\tau|^2 + 
\mathcal{O}\right) \right) \partial_r \gb . \]
This is easily evaluated from the expansion of $\gb$ recalled above and one gets
\begin{equation}\begin{split} \gb(\bar{\nabla}\nu(r),\cdot) \ = \ & \exp^{2r}\eta^2 + 
\frac{1}{2}\exp^r\gamma + \frac{R}{2}\exp^{-r} \gamma(J_0 \tau \cdot, \cdot) + 
\frac{1}{2} \exp^{-r}\gamma (\nabla_{\xi}\tau \cdot, \cdot ) \\
& + \frac{1}{3}\left(\frac{R^2}{8} - |\tau|^2 + \mathcal{O}\right)\,\eta^2  - 
\frac{1}{6}\left(\frac{R^2}{8} - |\tau|^2 + \mathcal{O}\right)\,\exp^{-r}\gamma .
\end{split}\end{equation}
One step further, this yields the (endomorphism) shape operator, which we shall develop
as 
\begin{equation}\label{II}
\II = \II_0 + \exp^{-r}\II_1 + \exp^{-\frac{3}{2}r}\II_{\frac{3}{2}} + 
\exp^{-2r}\II_2 + o(\exp^{-2r}), \end{equation} 
where $\II_0 = \Id_{\xi} + \frac{1}{2}\Id_H$, 
$\II_1 \ = \ \frac{R}{4}\Id_H \ - \ J_0\tau$, and
\[ \II_2 \ = \ \left(\frac{R^2}{8} - |\tau|^2 + \mathcal{O}\right)\,\Id_{\xi} \
+ \ \left(\frac{R^2}{16} + \frac{5}{2}|\tau|^2 + \mathcal{O}\right)\,\Id_{H} \ + 
\ \nabla_{\xi}\tau \]
and the precise value of $\II_{\frac{3}{2}}$ is irrelevant as before. As a {\sl last 
step}, we will obtain below
the desired contributions of the integral terms in Formula (\ref{chi}) by chasing the 
order $\exp^{-2r}$-terms. 

\medskip

A first easy consequence of the expression (\ref{II}) of $\II$ 
is that the $\nabla_{\xi}\tau$-term may
be forgotten: it would create a scalar term linear in $\nabla_{\xi}\tau$, and no
such scalar invariant exists.

\medskip

The contraction $\tv$ is now explicitly described as follows: if symmetric endomorphisms 
$A$, $B$, and $C$ are diagonal in a basis $\{e^0,e^1,e^2,e^3\}$ chosen to be
$\gb$-orthonormal and $J$-adapted, with eigenvalues $\alpha_r$, $\beta_s$ 
and $\gamma_t$, then:  
\begin{equation}
\tv(A\wedge B\wedge C) = \csum ( \alpha_r\beta_s\gamma_t )\,\varpi ,
\end{equation}
where $\csum$ denotes the sum over all permutations of $\{r,s,t\}$; moreover,
if $\rho$ is a curvature term (endomorphism on $2$-forms) with constant
coefficients in the same basis with diagonal entries $K_{rs} = 
<\rho(e^r\wedge e^s),e^r\wedge e^s>$, and $A$ is as above, then:
\begin{equation}
\tv(A\wedge\rho) = \asum (K_{rs}\lambda_{t})\,\varpi ,
\end{equation}
where $\asum$ denotes the sum over circular (not all) permutations of $\{r,s,t\}$. 

These formulae make easy the evaluation of all possible order $2$ terms but the one 
involving $W_2^-$: the computations are done in a basis $(\partial_r,\exp^{-r}\xi,
\exp^{-\frac{r}{2}}h,\exp^{-\frac{r}{2}}J_0h)$, with $h$ a $\gamma$-unit element of $H$
chosen to be an eigenvector of $J_0\tau$. This basis is orthogonal for $\gb$ except
for order $\frac{3}{2}$-terms, which we can neglect as usual, and for an order $2$-term 
due to the presence of $\nabla_{\xi}\tau$ in the expression of $\gb$. Taking into account
this last term yields a scalar correction linear in $\nabla_{\xi}\tau$, which must
necessarily vanish, hence one may also forget it. 
We summarize below the results, using the following notations: 
for any geometric quantity $s$, $s_k$ denotes the $k$-th order term in its asymptotic 
expansion; we also denote the order $2$ term in $\II$ as:
\[ \II_2 \ = \ A\, \Id_{\xi}\  + \ B\, \Id_H \ + \ \nabla_\xi\tau .\]
The term $\tv(\II\wedge\II\wedge\II)$ is the sum of contributions of type 
$<\varpi_a,\II_b\wedge\II_c\wedge\II_d>$, with $a+b+c+d=2$; the results are:
\begin{equation*}\begin{array}{|c|c|}
\hline
\textrm{\underline{involved terms}} & \textrm{\underline{result}} \\
\varpi_2, \II_0\wedge\II_0\wedge\II_0 & 
       1/16\, R^2 - 1/2\,|\tau|^2 \\
\varpi_1, \II_1\wedge\II_0\wedge\II_0 \ \textrm{(3 terms)} &
      - 3/4\, R^2 \\
\varpi_0, \II_1\wedge\II_1\wedge\II_0 \ \textrm{(3 terms)} &
       3/2\, R^2 - 6\, |\tau|^2\\
\varpi_0, \II_2\wedge\II_0\wedge\II_0 \ \textrm{(3 terms)} & 
      3/2\, A + 6 B \\
\hline
\end{array}\end{equation*}
And for the $\tv(\II\wedge\mathcal{R})$-term, the results are:
\begin{equation*}\begin{array}{|c|c|}
\hline
\textrm{\underline{involved terms}}  & \textrm{\underline{result}} \\
\varpi_2, \II_0\wedge\mathcal{R}_0 & 
    - 5/96\, R^2 + 5/12\, |\tau|^2 \\
\varpi_1, \II_1\wedge\mathcal{R}_0 & 
    1/16\, R^2\\
\varpi_0, \II_2\wedge\mathcal{R}_0 & - A - 1/2\, B \\
\hline
\end{array}\end{equation*}
For the last term, {\sl i.e.} $\tv(\II_0\wedge W_2^-)$, we have to rely on the following 
explicit expression of the second order correction to the curvature, extracted from 
\cite{ob-mh}. If $\omega_-^2 = 
e^0\wedge e^2 - e^1\wedge e^3$, and $\omega^3_- = e^0\wedge e^3 - e^1\wedge e^2$, then
\[ W_2^- = a\exp^{-2r} ((\omega^2_-)^2 - (\omega^3_-)^2) + b\exp^{-2r}
(\omega_-^2\omega^3_- + \omega^3_-\omega^2_-) ,\] 
$a$ and $b$ being reals.
The definition of $\tv$ then yields $\tv(\II_0\wedge W_2^-) = 0$. 
Putting together all the results obtained so far yields easily the expected Theorem.\qed

\bigskip

\section*{Appendix: variations of $R$ and $\tau$.}

\medskip

We give here a quick glimpse on the computations leading to the expression of the
variation of the curvature quantity $T(\eta) = \frac{R^2}{16} - \frac{5}{2}|\tau|^2$
of the Tanaka-Webster connection under a conformal deformation of the contact form
$\eta$. 

Let $\eta$ be a compatible contact form on $X$ and $f$ a positive function. We
denote by $(\ti,\tib)$ a local orthonormal (complex) coframe for $\eta$, {\it i.e.}
$d\eta = i \ti\wedge\tib$. The variation of the Tanaka-Webster curvature $R$ 
is well-known in dimension $3$ \cite{jerison-lee}; if $\widehat{R}$ is the curvature 
for $f^2\eta$, and $\Delta u = u_{,\bar1 1} + u_{,1 \bar 1}$ is the sub-elliptic 
Laplacian, then
\begin{equation}
\widehat{R} = f^{-3} \left( 2\,\Delta f + R\, f \right) .
\end{equation}
For the torsion $\tau$ and in lack of a precise reference, we will detail the 
computation a little bit. Starting from $(\ti,\tib)$, a local orthonormal coframe 
for $f^2\eta$ is then given by 
$\fti = f \left( \ti + 2i\, f_{,\bar1}\, \eta\right)$
and its complex conjugate. The Tanaka-Webster connection $1$-form $\omega_1^1$ and torsion
endomorphism $\tau_{\bar 1}^1$ for $\eta$ (resp. $f^2\eta$) are defined by
$d\ti \ = \ - \, \omega_1^1 \wedge \ti + \tau_{\bar 1}^1\eta\wedge\tib$
(resp. the same formula in the hatted version).
Computing at a point where $\omega_1^1$ is zero, one gets
\begin{equation}\begin{split}
d\fti & = df\wedge\ti + f\, d\ti + 2i f_{,\bar 1}\,d\eta + 
2i\, d\left( f_{,\bar 1}\right)\wedge\eta\\
      & = 3 f_{,\bar 1}\, \tib\wedge\ti + 2i\,( f_{,\bar1 1} + \frac{i}{2}f_{,0})
\,\ti\wedge\eta + ( 2if_{,\bar 1 \bar 1} - f \tau_{\bar 1}^{1})\, \tib\wedge\eta .
\end{split}\end{equation}
Identifying this with $-\, \hat{\omega}_1^1 \wedge \fti \ + 
\ \hat{\tau}_{\bar 1}^1\, \ftib \wedge f^2 \eta $, it comes finally:
\begin{equation}
\hat{\tau}_{\bar 1}^1 = f^{-2}\, \left( \tau_{\bar 1}^1 - 2i\, f^{-1} f_{,\bar 1\bar 1}
- 6i f^{-2} f_{,1}\, f_{,\bar 1}\right) .
\end{equation}
From these computations, the interested reader can easily derive the variation
of $T(\eta)$ under conformal changes in $\eta$.

\bigskip

\begin{small}
{\flushleft\sl Acknowledgements}. The author thanks Olivier Biquard, Gilles Carron and
Jean-Marc Schlenker for their interest in this work, and C. Robin Graham for useful 
comments.
\end{small}

\bigskip

\bibliographystyle{smfplain}

\medskip

\end{document}